\documentclass[fleqn,12pt]{article}

\usepackage{mathptmx}
\usepackage{mathrsfs}
\usepackage{amsmath}
\usepackage{amssymb}
\usepackage{amsthm}
\usepackage{bbm}
\usepackage{cite}
\usepackage{epstopdf}
\usepackage{geometry}
\usepackage{graphicx}
\usepackage{indentfirst}
\usepackage{multicol}
\usepackage{setspace}
\usepackage{subfigure}
\usepackage[T1]{fontenc}
\usepackage{lineno}
\usepackage{cite}
\usepackage{bm}
\usepackage[colorlinks,linkcolor=blue,anchorcolor=blue,citecolor=blue,CJKbookmarks=false]{hyperref}

\allowdisplaybreaks

\newtheorem{theorem}{Theorem}[section]
\newtheorem{remark}{Remark}[section]

\newtheorem{definition}{Definition}[section]

\begin{document}
\begin{sloppypar}
\begin{flushleft}
\LARGE{Optimal defined contribution pension management with jump diffusions and common shock dependence}
\end{flushleft}

\begin{flushleft}
{{Xiaoyi Zhang$^{1, }$}\renewcommand{\thefootnote}{\ddag}\footnote{E-mail: zhangxiaoyi19902@163.com.}
~|~
{Linlin Tian$^{2, }$}\renewcommand{\thefootnote}{\S}\footnote{Corresponding author. E-mail address: linlin.tian@dhu.edu.cn.}\\
\small  1. School of Economics and Management, Hebei University of Technology, Tianjin 300401, P.R. China.\\
\small  2. College of Science, Donghua University, Shanghai 201620, P.R. China.
}\,
\end{flushleft}

\noindent \textbf{Abstract}

\noindent This work deals with an optimal asset allocation problem for a defined contribution (DC) pension plan during its accumulation phase. The contribution rate is proportional to the individual's salary, the dynamics of which follows a Heston stochastic volatility model with jumps, and there are common shocks between the salary and the volatility. Since the time horizon of pension management might be long, the influence of inflation is considered in the context. The inflation index is subjected to a Poisson jump and a Brownian uncertainty. The pension plan aims to reduce fluctuations of terminal wealth by investing the fund in a financial market consisting of a riskless asset and a risky asset. The dynamics of the risky asset is given by a jump diffusion process. The closed form of investment decision is derived by the dynamic programming approach.

\noindent \textbf{Keywords}

\noindent DC pension plan, Stochastic volatility, Poisson process, Common shock dependence, Inflation, Hamilton-Jacobi-Bellman equation

\section{Introduction}

Pension fund is an important financial instrument for individuals to reallocate incomes and sustain consumption after retirement. Generally, according to determination of benefits, there are two typical kinds of pension plans: defined benefit (DB) and defined contribution (DC) pension plan. In DB plan, benefits are fixed in advance, while in DC case, contributions are fixed by the trustee. There are two phases for a pension scheme: the accumulation phase, which is the period from entry time to retirement time, and the decumulation phase, which is the period from retirement time to death time.

In the accumulation phase of DC pension scheme, the contributor contributes part of his or her salary to the fund. Since the salary is related to the profitability of the company, both works of Bodie et al.\cite{Bodie} and Dybvig and Liu \cite{Dybvig} assume that the salary process is spanned by the stock price. In addition, Guan and Liang \cite{Guan} and Li and Wang \cite{WangandLi} describe salary process by a Heston stochastic volatility model, i.e., salary is correlated with the volatility of the stock. Furthermore, Zeng et al.\cite{Zeng} assumes that the salary process is related to the stochastic volatility. It would be more realistic to introduce an independent random process on the stochastic salary process.

It is natural to insert a jump process in the stochastic salary due to promotion and job-hopping. Moreover, it is also realistic to introduce jumps in the volatility, which influenced by some unexpected events, such as economic crisis and policy adjustments by the government. In our model, the contribution rate of the pension scheme is proportional to the salary of the individual, the dynamics of which follows a Heston stochastic volatility model with jumps. In addition, salary and variance are correlated by means of a common shock. In reality, a common component may depict an event which has impact on both the salary and the volatility. Common shock models are widely used in insurance area. For instance, in Liang et al.\cite{LiangandYuenandZhang16}, the insurance risk model is modulated by a compound Poisson process, and the two jump-number processes are correlated through a common shock. Liang et al.\cite{LiangandYuenandZhang18} assumes that jumps in both the risky asset and insurance risk process are correlated through common shock dependence.

Since the period of a pension scheme is usually long, inflation risk has attracted increasing attention from academic aspect in the past decade. A number of studies concentrate on stochastic optimization problem for DC pension plan under inflation risk. For example, by considering inflation in the wealth process, Han and Hung \cite{Han} investigates the optimal asset allocation problem by dynamic programming approach. By assuming the price index follows a mean reversion model, Yao et al.\cite{Yao} solves an optimal portfolio decision problem using the mean-variance criterion. Other works on optimal control under inflation risk can be found in \cite{Chen} and \cite{Tang}.

As for the nominal price level of a representative bundle of commodity goods in the market, instead of using a continuous process, it would be realistic to introduce a jump diffusion model to reflect sudden shocks in the price index. Thus the dynamics of price index given by Zhang et al.\cite{Zhang07} and Zhang and Ewald \cite{Zhang10} is extended in our model, and a Poisson jump is included in the evolution of the index price.

For the evolution of stock price, jump diffusion model is widely used in asset allocation problems. Merton \cite{Merton} considers Poisson jumps in an optimal dynamic portfolio decision problem. In DC pension funding framework, Sun et al.\cite{Sun} deals with the precommitment and equilibrium investment strategies by incorporating jumps into the risky asset process. Works in DB pension management can be found in \cite{Delong}, \cite{XLiang} and \cite{Ngwira}. In this paper, we are interested in an optimal asset decision problem when the stock price is driven by a Brownian motion and a Poisson jump.

The rest of the paper is structured as follows. Section 2 describes the financial market with jump diffusion price index, as well as two tradable assets which are of interest for the pension management. This section also gives the pension model. The salary follows a Heston stochastic volatility model with jumps, and there exists common shock dependence between the salary and the volatility. Section 3 deals with a stochastic optimal control problem in order to minimize the fluctuation of the final real wealth over a finite horizon. Closed form investment strategy is given by solving the Hamilton-Jacobi-Bellman (HJB) equation. Finally, Section 4 establishes the conclusion.

\section{Model assumptions and notations}

Consider a probability space $(\Omega,\mathscr{F},\mathbb{P})$, with $\mathbb{P}$ the real world probability measure on $\Omega$ and $\mathscr{F}=\mathscr{F}^W\vee \mathscr{F}^{N}$. The filtration $\mathscr{F}^W=\left\{\mathscr{F}^W_t\right\}_{t \geq 0}$ is generated by a five dimensional standard Brownian motion $(W_r,W_S,W_L,W_V,W_\Pi)$, i.e., $\mathscr{F}^W_t=\sigma\left\{(W_r(s),W_S(s),W_L(s),W_V(s),W_\Pi(s));0\leq s\leq t\right\}, t\geq0$. $W_r$ and $W_\Pi$ are correlated, which is captured by the coefficient $\rho_{\Pi r}\in(-1,1)$. $W_L$ and $W_V$ are also correlated, which is captured by the coefficient $\rho_{LV}\in(-1,1)$. Let the filtration $\mathscr{F}^N=\left\{\mathscr{F}^N_t\right\}_{t \geq 0}$ be generated by a five dimensional Poisson process $(N_S,N_L,N_c,N_V,N_\Pi)$ with intensity $(\lambda_S,\lambda_L,\lambda_c,\lambda_V,\lambda_\Pi)$, where $\lambda_S,\lambda_L,\lambda_c,\lambda_V,\lambda_\Pi\in \mathbb{R}^+$, i.e., $\mathscr{F}^N_t=\sigma\left\{N_S(s),N_L(s),N_c(s),N_V(s),N_\Pi(s);0\leq s\leq t\right\}, t\geq0$. Suppose that Poisson processes are mutually independent. Besides, Brownian motions are independent of Poisson processes on $(\Omega,\mathscr{F},\mathbb{P})$.

\subsection{The financial market}

Assume that the market consists of two underlying instruments which are traded continuously over time and perfectly divisible. Following the work of Eisenberg\cite{Eisenberg}, we assume that the discount factor is a geometric Brownian motion:
\begin{eqnarray}\label{1}
\exp\big[r+mt+\zeta W_r(t)\big],
\end{eqnarray}
where $r, m \in \mathbb{R}^+$, and $\zeta \in \mathbb{R}$. Thus the riskless asset $S_0(t)$ evolves according to the the following form
\begin{eqnarray}
\begin{array}{rcl}\label{2}
$$\dfrac{dS_0(t)}{S_0(t)}=(m+\dfrac{\zeta^2}{2})dt+\zeta dW_r(t),$$
\end{array}
\end{eqnarray}
with initial price $S_0(0)=e^r$.

\begin{remark}
In most of previous works, it is assumed that the cash account has the following form
\begin{eqnarray}\label{3}
\dfrac{dS_0(t)}{S_0(t)}=r(t)dt,
\end{eqnarray}
with initial price $S_0(0)=1$, and the interest rate $r(t)$ follows the CIR or Vasicek model, see, for example, Guan and Liang\cite{Guan}. However, in our model, Eq.(\ref{2}) allows us to find the explicit solution, and the randomness of the interest rate is remained.
\end{remark}

Besides the cash account, the trustee also has the opportunity to invest the fund into a stock with the dynamics
\begin{eqnarray}\label{4}
\dfrac{dS(t)}{S(t-)}=\mu_S(t)dt+\sigma_{SS}(t)dW_S(t)+\eta_S(t)dN_S(t),
\end{eqnarray}
where $\mu_S(t)$ is the appreciation rate for the stock. $\sigma_{SS}(t)$ is the volatility associated with the diffusion component of the stock price. $\eta_S(t)$ denotes the magnitude of a jump. We state that $\eta_S(t)>-1$ to prevent the process jumping to a value below zero. $W_S$ describes the fluctuation, and $N_S$ describes the jump of the stock price. For simplicity, it is assumed that $W_S$ and $N_S$ are independent stochastic processes.

\subsection{The pension model}

This paper considers the accumulation phase of a DC type pension plan. Assume that the entry time of a pensioner is the initial time $0$, and his/her retirement time is the terminal time in our model. Denote the pensioner's death time as $\tau$, which is a positive random variable defined on the probability space $(\Omega,\mathscr{F},\mathbb{P})$. The mortality rate $\lambda(t)$ is defined as
\begin{eqnarray}\label{5}
\lambda(t)=\lim \limits_{\triangle t \rightarrow 0 } \dfrac{ \mathbb{P}(t < \tau < t+\triangle t\;|\;\tau > t)}{\triangle t}.
\end{eqnarray}

In a general pension plan, the pensioner pays contributions before the retirement time $T$, where $T\in R^+$. The level of the contribution rate is usually defined as a proportion $\xi(0\leq\xi\leq1)$ of the pensioner's salary. In previous works, for instance, \cite{Guan}, \cite{WangandLi} and \cite{Zeng}, it is assumed that the stochastic salary is driven by the Heston's SV model, i.e., the salary process has stochastic volatility, and the salary return variance is governed by a mean-reverting process. However, the salary is only influenced by the interest rate, the stock price and the volatility, which is not quite realistic. Here we add a Brownian motion $W_L$ in the salary process to describe the fluctuation of the salary itself. In addition, we assume that both the salary process and the stochastic variance have jumps.

The dynamics of salary $L(t)$ is given by the following differential equation
\begin{eqnarray}\label{6}
\dfrac{dL(t)}{L(t-)}=\mu_L(t)dt+\sigma_{LS}(t)dW_S(t)+\sqrt{V(t)}dW_L(t)+\eta_{LL}(t)dN_L(t)+\eta_{Lc}(t)dN_c(t),
\end{eqnarray}
with initial value $L(0)=L_0$. $\mu_L(t)$ denotes the instantaneous expected rate of salary. $\sigma_{LS}(t)$ is the instantaneous volatility scale factor measuring how risk source of stock price affect the salary. $V(t)$ is the stochastic volatility which is explained in Eq.(\ref{7}). Here we assume that there exists common shock dependence between the salary and the variance. $\eta_{LL}(t)(>-1)$ and $\eta_{Lc}(t)(>-1)$ denote magnitude of jumps associated with Poisson processes $N_L(t)$ and $N_c(t)$, respectively. $N_L$ describes the jump of the salary itself, and $N_c$ describes the common shock between the salary (given by Eq.(\ref{6})) and the stochastic volatility $V(t)$ (given by Eq.(\ref{7})). In Eq.(\ref{6}), Browian motions and Poisson processes are mutually independent.

The dynamics of the stochastic volatility $V(t)$ is given by
\begin{eqnarray}\label{7}
dV(t)=\kappa(\delta-V(t))dt+\sigma_V\sqrt{V(t)}dW_V(t)+\eta_{VV}(t)dN_V(t)+\eta_{Vc}(t)dN_c(t),
\end{eqnarray}
where $\kappa$ denotes the mean-reversion rate, $\delta$ denotes the long-run mean, and $\sigma_V$ is the volatility coefficient. There is a condition $2\kappa\delta>\sigma_V^2$ to guarantee the volatility process $V(t)>0$. The Brownian motion $W_V$ describes the fluctuation of the volatility, and the Poisson process $N_V$ describes the jump of the volatility. We state that $\eta_{VV}>-1$ and $\eta_{Vc}>-1$ to prevent the process $V(t)$ jumping to a value below zero.

\section{The optimal portfolio}

The aim of the stochastic control problem is to find the optimal investment decision. The pension trustee continuously decides weights invested into the cash account and the stock. Denote the nominal wealth at time $t$ as $X(t)$. Under the investment policy chosen, it is easy to get the following stochastic differential equation which describes the evolution of the wealth
\begin{equation}\label{8}
dX(t)=X(t)(1-\pi(t))\dfrac{dS_0(t)}{S_0(t)}+X(t-)\pi(t)\dfrac{dS(t)}{S(t-)}+\xi(t)L(t)dt,
\end{equation}
with $X(0)=X_0>0$. $\pi(t)$ denotes the weight invested into the stock at time $t$. The remainder, $1-\pi(t)$ is the proportion invested into the cash account. Borrowing and short-selling are permitted in the context. A negative value of $\pi(t)$ means that the pension trustee takes a short position in the stock, while a negative value of $1-\pi(t)$ reflects that the trustee borrows money from the bank to purchase the risky asset.

By substituting Eq.(\ref{2}) and Eq.(\ref{4}) into Eq.(\ref{8}), we obtain that
\begin{equation}
\begin{array}{rcl}\label{9}
$$dX(t)\!\!\!\!&=&\!\!\!\!X(t)\Big[(m+\dfrac{\zeta^2}{2})+\pi(t)(\mu_S(t)-(m+\dfrac{\zeta^2}{2}))\Big]dt+\xi(t)L(t)dt\\

\!\!\!\!&{}&\!\!\!\!+X(t)(1-\pi(t))\zeta(t)dW_r(t)+X(t)\pi(t)\sigma_{SS}(t)dW_S(t)+X(t-)\pi(t-)\eta_S(t)dN_S(t).$$
\end{array}
\end{equation}

As mentioned in Section 1, the time horizon for the accumulation phase of a pension fund (from time $0$ to $T$) is usually long, hence the influence of inflation is necessary to be considered in the context.

The price index at time $t$ is denoted by $\Pi(t)$ and the dynamics is driven by a jump diffusion process of the following type
\begin{eqnarray}\label{10}
\dfrac{d\Pi(t)}{\Pi(t-)}=\mu_\Pi(t)dt+\sigma_\Pi(t)dW_\Pi(t)+\eta_\Pi(t)dN_\Pi(t),
\end{eqnarray}
with initial value $\Pi(0)=\Pi_0 > 0$. $\mu_\Pi(t)$ is the instantaneous expected inflation rate. $\sigma_\Pi(t)$ is the instantaneous volatility associated with the diffusion component and $\eta_\Pi(t)$ denotes the magnitude of a jump with the condition $\eta_\Pi > -1$ to ensure that the price index remains strictly positive. In Eq.(\ref{10}), $W_\Pi$ and $N_\Pi$ are independent stochastic processes.

Next we define the corresponding real salary process as the following:

\begin{definition}
The real salary process is defined by
\begin{equation}\label{11}
\overline{L}(t)=\dfrac{L}{\Pi}(t).
\end{equation}
\end{definition}

Applying the quotient rule of the It\^{o}'s formula, $\overline{L}$ is given by
\begin{equation}
\begin{array}{rcl}\label{12}
$$d\overline{L}(t)\!\!\!\!&=&\!\!\!\!d\Big[\dfrac{L}{\Pi}\Big](t)\\

\!\!\!\!&=&\!\!\!\!\overline{L}(t)(\mu_L(t)-\mu_\Pi(t)+\sigma_\Pi^2(t))dt
+\overline{L}(t)\sigma_{LS}(t)dW_S(t)+\overline{L}(t)\sqrt{V(t)}dW_L(t)\\

\!\!\!\!&{}&\!\!\!\!-\overline{L}(t)\sigma_\Pi(t)dW_\Pi(t)
+\overline{L}(t-)\eta_{LL}(t)dN_L(t)+\overline{L}(t-)\eta_{Lc}(t)dN_c(t)\\

\!\!\!\!&{}&\!\!\!\!+\overline{L}(t-)(\eta_\Pi^2(t)-\eta_\Pi(t))dN_\Pi(t).$$
\end{array}
\end{equation}
with initial value $\overline{L}(0)=L_0/\Pi_0\triangleq\overline{L}_0$.

Then the real wealth process with the consideration of inflation follows
\begin{equation}
\begin{array}{rcl}\label{13}
$$d\overline{X}(t)\!\!\!\!&=&\!\!\!\!d\Big[\dfrac{X}{\Pi}\Big](t)\\

$$\!\!\!\!&=&\!\!\!\!\overline{X}(t)\Big[(m+\dfrac{\zeta^2}{2})+\pi(t)(\mu_S(t)-(m+\dfrac{\zeta^2}{2})
+\zeta\sigma_\Pi(t)\rho_{\Pi r}(t))-\mu_\Pi(t)+\sigma_\Pi^2(t)\\

\!\!\!\!&{}&\!\!\!\!-\zeta\sigma_\Pi(t)\rho_{\Pi r}(t)\Big]dt+\xi(t)\overline{L}(t)dt
+\overline{X}(t)(1-\pi(t))\zeta dW_r(t)+\overline{X}(t)\pi(t)\sigma_{SS}(t)dW_S(t)\\

\!\!\!\!&{}&\!\!\!\!-\overline{X}(t)\sigma_\Pi(t)dW_\Pi(t)
+\overline{X}(t-)\pi(t)\eta_S(t)dN_S(t)+\overline{X}(t-)(\eta_\Pi^2(t)-\eta_\Pi(t))dN_\Pi(t),
\end{array}
\end{equation}
with initial condition $\overline{X}(0)=X_0/\Pi_0\triangleq\overline{X}_0$.

Next we restrict the strategies in order to fulfil some technical conditions. A strategy $\pi(\cdot)$ is a control process which is $\mathscr{F}_t$-measurable, Markovian and stationary, satisfying the condition
\begin{equation}\label{14}
\mathbb{E}\left\{\displaystyle{\int_0^\infty \pi^2(t)dt} \right\} < \infty.
\end{equation}

Denote $\mathscr{A}_{\overline{X}_0,\overline{L}_0}$ the set of all admissible controls, i.e., it is the set of all measurable processes $\left\{{\pi}(t)\right\}_{t\geq 0}$, which satisfies Eq.(\ref{14}).

Assume that the pension trustee has a preference to minimize the expected value of the fluctuations of the terminal wealth until time $\tau \wedge T$, where $T$ is the terminal time of the control problem. The objective is to minimize
\begin{equation}\label{15}
J(t,\overline{X},\overline{L},V)
=\mathbb{E}_t\bigg[\big[\alpha_1+\beta_1(\overline{X}(T)-X_1^*)\big]^2\cdot\mathbbm{1}_{\left\{\tau>T\right\}}
+\big[\alpha_2+\beta_2(\overline{X}(\tau)-X_2^*)\big]^2\cdot\mathbbm{1}_{\left\{\tau\leq T\right\}} \Big| \tau>t \bigg],
\end{equation}
with $\mathbb{E}_t$ the conditional expectation given the filtration $\left\{\mathscr{F}_t\right\}_{t \geq 0}$. $X_1^*$ and $X_2^*$ are target funds of the plan at time $T$ and $\tau$, respectively. The deviation between the actual fund and the target fund is called the discontinuity risk, see Wang et al.\cite{Wang}.

In quadratic loss functions, any deviations from $\overline{X}$ are penalized when the wealth is different from the target $X_1^*$ or $X_2^*$, so that a cost measured by the loss function must be paid. It is assumed that $\alpha_1,\alpha_2>0$ and $\beta_1,\beta_2<0$ in Eq.(\ref{15}), thus, under funding is more penalized than over funding, see, for example, Devolder, Janssen, and Manca\cite{Devolder} and Zhang and Guo\cite{ZhangangGuo}.

According to
\begin{equation}\label{16}
\mathbb{E}_t\bigg[\big[\alpha_1+\beta_1(\overline{X}(T)-X_1^*)\big]^2\cdot\mathbbm{1}_{\left\{\tau>T\right\}}\Big| \tau>t\bigg]
=\mathbb{E}_t\bigg[\big[\alpha_1+\beta_1(\overline{X}(T)-X_1^*)\big]^2 e^{-\int_t^T\lambda(u)du}\bigg],
\end{equation}
and
\begin{equation}\label{17}
\mathbb{E}_t\bigg[\big[\alpha_2+\beta_2(\overline{X}(\tau)-X_2^*)\big]^2\cdot\mathbbm{1}_{\left\{\tau\leq T\right\}}\Big| \tau>t\bigg]
=\mathbb{E}_t\bigg[\displaystyle{\int_t^T}\big[\alpha_2+\beta_2(\overline{X}(s)-X_2^*)\big]^2\lambda(s) e^{-\int_t^s\lambda(u)du}ds\bigg],
\end{equation}
the objective function with an uncertain lifetime is converted into the following deterministic horizonal function
\begin{equation}
\begin{array}{rcl}\label{18}
$$\!\!\!\!&{}&\!\!\!\!J(t,\overline{X},\overline{L},V)\\

\!\!\!\!&=&\!\!\!\!\mathbb{E}_t\bigg[\displaystyle{\int_t^T}
\big[\alpha_2+\beta_2(\overline{X}(s)-X_2^*)\big]^2\lambda(s)e^{-\int_t^s\lambda(u)du}ds
+\big[\alpha_1+\beta_1(\overline{X}(T)-X_1^*)\big]^2 e^{-\int_t^T\lambda(u)du}\bigg].$$
\end{array}
\end{equation}

The dynamic programming approach is used to solve the stochastic optimization problem. Define the value function as
\begin{equation}\label{19}
\varphi(t,\overline{X},\overline{L},V)=\min \limits_{\left\{\pi\right\}}\left\{ J(t,(\overline{X},\overline{L},V);\pi): \mbox{subject\;to}\;(\ref{13}), (\ref{12}), (\ref{7})\right\}.
\end{equation}

In stochastic optimal control theory, the HJB equation accomplishes the connection between the value function and the optimal control, see, for instance, Fleming and Soner\cite{Fleming} and Yong and Zhou \cite{Yong}:
\begin{equation}\label{20}
\min\limits_{\left\{\pi\right\}}\Psi(\pi)=0,
\end{equation}
where
\begin{equation}
\begin{array}{rcl}\label{21}
$$\!\!\!\!&{}&\!\!\!\!\Psi(\pi)\\

\!\!\!\!&=&\!\!\!\!\varphi_t+\lambda\big[\alpha_2+\beta_2(\overline{X}-X_2^*)\big]^2
-\lambda \varphi+\varphi_{\overline{X}}\overline{X}\Big[(m+\dfrac{\zeta^2}{2})
+\pi(\mu_S-(m+\dfrac{\zeta^2}{2})+\zeta\sigma_\Pi\rho_{\Pi r})\\

\!\!\!\!&{}&\!\!\!\!-\mu_\Pi+\sigma_\Pi^2-\zeta\sigma_\Pi\rho_{\Pi r}\Big]+\varphi_{\overline{X}}\xi\overline{L}
+\dfrac{1}{2}\varphi_{\overline{X}\;\overline{X}}\overline{X}^2\Big[(1-\pi)^2\zeta^2+\pi^2\sigma_{SS}^2+\sigma_\Pi^2-2(1-\pi)\\

\!\!\!\!&{}&\!\!\!\!\cdot\zeta\sigma_\Pi\rho_{\Pi r}\Big]+\varphi_{\overline{L}}\overline{L}(\mu_L-\mu_\Pi+\sigma_\Pi^2)
+\dfrac{1}{2}\varphi_{\overline{L}\;\overline{L}}\overline{L}^2 (\sigma_{LS}^2+V+\sigma_\Pi^2)
+\varphi_V\kappa(\delta-V)+\dfrac{1}{2}\varphi_{VV}\sigma_V^2V\\

\!\!\!\!&{}&\!\!\!\!+\varphi_{\overline{X}\;\overline{L}}\overline{X}\;\overline{L}\Big[\pi\sigma_{SS}\sigma_{LS}
+\sigma_\Pi^2-(1-\pi)\zeta\sigma_\Pi\rho_{\Pi r}\Big]\!+\!\varphi_{\overline{L}V}\overline{L}V\sigma_V\rho_{LV}
\!+\!\lambda_S\Big[\varphi(t,\overline{X}(1\!+\!\pi\eta_S),\overline{L},V)\\

\!\!\!\!&{}&\!\!\!\!-\varphi(t,\overline{X},\overline{L},V)\Big]
+\lambda_L\Big[\varphi(t,\overline{X},\overline{L}(1+\eta_{LL}),V)-\varphi(t,\overline{X},\overline{L},V)\Big]
+\lambda_V\Big[\varphi(t,\overline{X},\overline{L},V+\eta_{VV})\\

\!\!\!\!&{}&\!\!\!\!-\varphi(t,\overline{X},\overline{L},V)\Big]
+\lambda_C\Big[\varphi(t,\overline{X},\overline{L}(1+\eta_{Lc}),V+\eta_{Vc})-\varphi(t,\overline{X},\overline{L},V)\Big]\\

\!\!\!\!&{}&\!\!\!\!+\lambda_\Pi\Big[\varphi(t,\overline{X}(1+(\eta_\Pi^2-\eta_\Pi)),\overline{L}(1+(\eta_\Pi^2-\eta_\Pi)),V)
-\varphi(t,\overline{X},\overline{L},V)\Big],$$
\end{array}
\end{equation}
with terminal condition $\varphi(T,\overline{X},\overline{L},V)=\big[\alpha_1+\beta_1(\overline{X}(T)-X_1^*)\big]^2$. $\varphi_t$, $\varphi_{\overline{X}}$, $\varphi_{\overline{L}}$, $\varphi_{V}$, $\varphi_{\overline{X}\;\overline{X}}$, $\varphi_{\overline{L}\;\overline{L}}$, $\varphi_{VV}$, $\varphi_{\overline{X}\;\overline{L}}$ and $\varphi_{\overline{L}V}$denote the first and second order partial derivatives of the value function $\varphi$ with respect to $t$, $\overline{X}$, $\overline{L}$ and $V$, respectively.

If there exist a twice continuously differentiable solution of Eq.(\ref{21}), strictly convex, then the minimizer of the investment strategy is obtained by the optimal functional $\pi^{\ast}$, which satisfies the following necessary conditions
\begin{equation}\label{22}
\Psi(\pi^{\ast})=0,
\end{equation}
\begin{equation}\label{23}
\dfrac{d\Psi}{d\pi}(\pi^{\ast})=0.
\end{equation}

We shall frequently use the following notations. Define
\begin{equation}\label{24}
\varpi_1=\mu_S-(m+\dfrac{\zeta^2}{2})+\zeta\sigma_\Pi\rho_{\Pi r}+\lambda_S\eta_S,
\end{equation}
\begin{equation}\label{25}
\varpi_2=\zeta\sigma_\Pi\rho_{\Pi r}+\sigma_{SS}\sigma_{LS},
\end{equation}
\begin{equation}\label{26}
\varpi_3=\zeta\sigma_\Pi\rho_{\Pi r}-\zeta^2,
\end{equation}
\begin{equation}\label{27}
\varpi_4=\zeta^2+\sigma_{SS}^2+\lambda_S\eta_S^2.
\end{equation}

By using the first order condition and solving the HJB equation, the explicit form of the optimal investment decision is given by the following theorem.
\begin{theorem}\label{th3.1}
(Main result) The optimal investment strategy on the stock is given by
\begin{equation}\label{28}
\pi^{\ast}(t) =-\dfrac{2\varphi_1(t)\overline{X}+\varphi_2(t)+\varphi_5(t)\overline{L}}{2\varphi_1(t)\overline{X}}\cdot\dfrac{\varpi_1}{\varpi_4}
-\dfrac{\varphi_5(t)\overline{L}}{2\varphi_1(t)\overline{X}}\cdot\dfrac{\varpi_2}{\varpi_4}-\dfrac{\varpi_3}{\varpi_4}.
\end{equation}

The value function is
\begin{equation}\label{29}
\varphi(t,\overline{X},\overline{L},V)=\varphi_1(t)\overline{X}^2+\varphi_2(t)\overline{X}+\varphi_3(t,V)\overline{L}^2
+\varphi_4(t)\overline{L}+\varphi_5(t)\overline{X}\;\overline{L}+\varphi_6(t).
\end{equation}

In above equations,
\begin{equation}
\begin{array}{rcl}\label{30}
$$\varphi_1(t)=\lambda\beta_2^2\displaystyle\int_t^Te^{\int_t^s a_1(u)du}ds
+\beta_1^2e^{\int_t^T a_1(s)ds},$$
\end{array}
\end{equation}
\begin{equation}
\begin{array}{rcl}\label{31}
$$\varphi_2(t)=2\lambda(\alpha_2\beta_2-\beta_2^2X_2^*)\displaystyle\int_t^Te^{\int_t^s a_2(u)du}ds
+2\beta_1(\alpha_1-\beta_1 X^*_1)e^{\int_t^T a_2(s)ds},$$
\end{array}
\end{equation}
\begin{equation}
\begin{array}{rcl}\label{32}
$$\varphi_3(t,V)=\displaystyle{\int_t^T}\widetilde{\varphi}_{31}(t;\tau)e^{\widetilde{\varphi}_{32}(t;\tau)V}d\tau,$$
\end{array}
\end{equation}
\begin{equation}
\begin{array}{rcl}\label{33}
$$\varphi_4(t)=\displaystyle{\int_t^T} e^{\int_0^s [a_4(u)+\lambda_c\eta_{Lc}]du}\Big[\xi\varphi_2(s)
-\dfrac{\varphi_2(s)\varphi_5(s)}{2\varphi_1(s)}\cdot\dfrac{\varpi_1(\varpi_1+\varpi_2)}{\varpi_4}\Big]ds\cdot e^{-\int{^t_0}[a_4(s)+\lambda_c\eta_{Lc}]ds},$$
\end{array}
\end{equation}
\begin{equation}
\begin{array}{rcl}\label{34}
$$\varphi_5(t)=2\displaystyle{\int_t^T} e^{\int_0^s [a_5(u)+\lambda_c\eta_{Lc}]du}\xi\varphi_1(s)ds\cdot e^{-\int{^t_0}[a_5(s)+\lambda_c\eta_{Lc}]ds},$$
\end{array}
\end{equation}
\begin{equation}
\begin{array}{rcl}\label{35}
$$\varphi_6(t)=\displaystyle\int_t^T e^{\int_t^s a_1(u)du}\Big[\lambda(\alpha_2-\beta_2 X^*_2)^2
-\dfrac{\varphi_2^2}{4\varphi_1}\cdot\dfrac{\varpi_1^2}{\varpi_4}\Big]ds+(\alpha_1-\beta_1X^*_1)^2e^{\int_t^T a_1(s)ds},$$
\end{array}
\end{equation}
where $a_1$, $a_2$, $\widetilde{\varphi}_{31}$, $\widetilde{\varphi}_{32}$, $a_4$ and $a_5$ are given by Eq.(\ref{41}), Eq.(\ref{42}), Eq.(\ref{80}), Eq.(\ref{77}), Eq.(\ref{44}) and Eq.(\ref{45}), respectively.
\end{theorem}

\emph{Proof.} See Appendix.  $\square$

\section{Conclusion}

This paper analyses the optimal investment strategy for a DC type pension scheme during its accumulation phase, where the price of the risky asset follows a jump diffusion process. As the management phase is usually long, inflation is considered in the context, and the price index follows a jump diffusion process. We assume that dynamics of salary follows a Heston stochastic volatility model with jumps, and there are common shocks between the salary and the volatility. In order to determine the investment strategy and minimize the fluctuations of terminal wealth, the dynamic programming technique is applied in our work. The explicit solution of the optimal problem is derived from the correlated HJB equation.

Now we point out some further research. The first is to consider the risky asset follows a Heston stochastic volatility model with jumps, and there are common shocks between the stock and the volatility. As the stock market and salary process may have similar jumps caused by the same event, the second is to consider common shocks between the salary and the stock price.\\

\noindent \textbf{Appendix}

The proof of Theorem \ref{th3.1}.

From the HJB equation (\ref{20}), we conjecture that its solution takes a quadratic homogeneous form with $\varphi\in C^{1,2}$ and $\varphi_{\overline{X}\;\overline{X}}>0$, as the following
\begin{equation}\label{36}
\varphi(t,\overline{X},\overline{L},V)=\varphi_1(t,V)\overline{X}^2+\varphi_2(t,V)\overline{X}+\varphi_3(t,V)\overline{L}^2
+\varphi_4(t,V)\overline{L}+\varphi_5(t,V)\overline{X}\;\overline{L}+\varphi_6(t,V),
\end{equation}
where $\varphi_1(\cdot,\cdot)$, $\varphi_2(\cdot,\cdot)$, $\varphi_3(\cdot,\cdot)$, $\varphi_4(\cdot,\cdot)$, $\varphi_5(\cdot,\cdot)$, and $\varphi_6(\cdot,\cdot)$ are six suitable functions with terminal conditions $\varphi_1(T,V)=\beta_1^2$, $\varphi_2(T,V)=2\beta_1(\alpha_1-\beta_1 X^*_1)$, $\varphi_3(T,V)=\varphi_4(T,V)=\varphi_5(T,V)=0$ and $\varphi_6(T,V)=(\alpha_1-\beta_1X^*_1)^2$.

Differentiating Eq.(\ref{36}) with respect to $t$, $\overline{X}$, $\overline{L}$, $V$, we have
\begin{equation}
\begin{array}{rcl}\label{37}
&&\varphi_t=\varphi_{1t}\overline{X}^2+\varphi_{2t}\overline{X}+\varphi_{3t}\overline{L}^2
+\varphi_{4t}\overline{L}+\varphi_{5t}\overline{X}\;\overline{L}+\varphi_{6t},\qquad
\varphi_{\overline{X}}=2\varphi_1 \overline{X}+\varphi_2+\varphi_5\;\overline{L},\\

&&\varphi_{\overline{X}\;\overline{X}}=2\varphi_1,\qquad
\varphi_{\overline{L}}=2\varphi_3 \overline{L}+\varphi_4+\varphi_5\;\overline{X},\qquad \varphi_{\overline{L}\;\overline{L}}=2\varphi_3,\qquad
\varphi_{\overline{X}\;\overline{L}}=\varphi_5,\\

&&\varphi_V=\varphi_{1V}\overline{X}^2+\varphi_{2V}\overline{X}+\varphi_{3V}\overline{L}^2
+\varphi_{4V}\overline{L}+\varphi_{5V}\overline{X}\;\overline{L}+\varphi_{6V},\\

&&\varphi_{VV}=\varphi_{1VV}\overline{X}^2+\varphi_{2VV}\overline{X}+\varphi_{3VV}\overline{L}^2
+\varphi_{4VV}\overline{L}+\varphi_{5VV}\overline{X}\;\overline{L}+\varphi_{6VV},
\end{array}
\end{equation}
where $\varphi_{1t}$, $\varphi_{1V}$ and $\varphi_{1VV}$ denote the first and second order derivatives of $\varphi_{1}$ with respect to $t$ and $V$, respectively. The derivatives of $\varphi_{2}$, $\varphi_{3}$, $\varphi_{4}$, $\varphi_{5}$ and $\varphi_{6}$ are defined in the same way.

Substituting Eq.(\ref{36}) and Eq.(\ref{37}) into Eq.(\ref{20}) and Eq.(\ref{21}), and rearranging terms by the order of $\pi$, we obtain that
\begin{equation}
\begin{array}{rcl}\label{38}
$$&&\!\!\!\!\!\!\!\!\!\!\min\limits_{\left\{\pi\right\}}
\varphi_1\overline{X}^2\varpi_4\pi^2+(2\varphi_1\overline{X}+\varphi_2+\varphi_5\;\overline{L})\overline{X}\varpi_1\pi
+\varphi_5\overline{X}\;\overline{L}\varpi_2\pi+2\varphi_1\overline{X}^2\varpi_3\pi+\varphi_1\overline{X}^2(\zeta^2\\

&&\!\!-2\zeta\sigma_\Pi\rho_{\Pi r})-\varphi_5\overline{X}\;\overline{L}\zeta\sigma_\Pi\rho_{\Pi r}
+\varphi_{1t}\overline{X}^2+\varphi_{2t}\overline{X}+\varphi_{3t}\overline{L}^2
+\varphi_{4t}\overline{L}+\varphi_{5t}\overline{X}\;\overline{L}+\varphi_{6t}+\lambda\beta_2^2\overline{X}^2\\

&&\!\!+2\lambda(\alpha_2\beta_2-\beta_2^2 X^*_2)\overline{X}+\lambda(\alpha_2-\beta_2 X^*_2)^2
-\lambda (\varphi_1\overline{X}^2+\varphi_2\overline{X}+\varphi_3\overline{L}^2
+\varphi_4\overline{L}+\varphi_5\overline{X}\;\overline{L}+\varphi_6)\\

&&\!\!+(2\varphi_1\overline{X}+\varphi_2+\varphi_5\;\overline{L})\overline{X}
\Big[(m+\dfrac{\zeta^2}{2})-\mu_\Pi+\sigma_\Pi^2-\zeta\sigma_\Pi\rho_{\Pi r}\Big]
+(2\varphi_1\overline{X}+\varphi_2+\varphi_5\;\overline{L})\overline{L}\xi\\

&&\!\!+\varphi_1\overline{X}^2\sigma_\Pi^2+(2\varphi_3\overline{L}+\varphi_4+\varphi_5\overline{X})\overline{L}
(\mu_L-\mu_\Pi+\sigma_\Pi^2)+\varphi_3\overline{L}^2(\sigma_{LS}^2+V+\sigma_\Pi^2)+(\varphi_{1V}\overline{X}^2\\

&&\!\!+\varphi_{2V}\overline{X}+\varphi_{3V}\overline{L}^2
+\varphi_{4V}\overline{L}+\varphi_{5V}\overline{X}\;\overline{L}+\varphi_{6V})\kappa(\delta-V)
+\dfrac{1}{2}(\varphi_{1VV}\overline{X}^2+\varphi_{2VV}\overline{X}+\varphi_{3VV}\overline{L}^2\\

&&\!\!+\varphi_{4VV}\overline{L}+\varphi_{5VV}\overline{X}\;\overline{L}+\varphi_{6VV})\sigma_V^2V
+\varphi_5\overline{X}\;\overline{L}\sigma_\Pi^2
+(2\varphi_{3V}\overline{L}+\varphi_{4V}+\varphi_{5V}\;\overline{X})\overline{L}V\sigma_V\rho_{LV}\\

&&\!\!+\varphi_3\overline{L}^2\lambda_L\eta_L^2
+(2\varphi_3 \overline{L}^2+\varphi_4\overline{L}+\varphi_5\;\overline{X}\;\overline{L} )\lambda_L\eta_L
+\lambda_V\Big[(\varphi_1(t,V+\eta_{VV})-\varphi_1(t,V))\overline{X}^2\\

&&\!\!+(\varphi_2(t,V+\eta_{VV})-\varphi_2(t,V))\overline{X}+(\varphi_3(t,V+\eta_{VV})
-\varphi_3(t,V))\overline{L}^2+(\varphi_4(t,V+\eta_{VV})\\

&&\!\!-\varphi_4(t,V))\overline{L}+(\varphi_5(t,V+\eta_{VV})-\varphi_5(t,V))\overline{X}\;\overline{L}
+(\varphi_6(t,V+\eta_{VV})-\varphi_6(t,V))\Big]\\

&&\!\!+\lambda_c\Big[(\varphi_1(t,V+\eta_{Vc})-\varphi_1(t,V))\overline{X}^2+(\varphi_2(t,V+\eta_{Vc})-\varphi_2(t,V))\overline{X}
+(\varphi_3(t,V+\eta_{Vc})\\

&&\!\!-\varphi_3(t,V))\overline{L}^2+\varphi_3(t,V+\eta_{Vc})(\eta_{Lc}^2+2\eta_{Lc})\overline{L}^2
+(\varphi_4(t,V+\eta_{Vc})-\varphi_4(t,V))\overline{L}\\

&&\!\!+\varphi_4(t,V+\eta_{Vc})\eta_{Lc}\overline{L}
+(\varphi_5(t,V+\eta_{Vc})-\varphi_5(t,V))\overline{X}\;\overline{L}+\varphi_5(t,V+\eta_{Vc})\eta_{Lc}\overline{X}\;\overline{L}\\

&&\!\!+(\varphi_6(t,V+\eta_{Vc})-\varphi_6(t,V))\Big]+\lambda_\Pi\Big[\varphi_1(t,V)((\eta_\Pi^2-\eta_\Pi)^2
+2(\eta_\Pi^2-\eta_\Pi))\overline{X}^2\\

&&\!\!+\varphi_2(t,V)(\eta_\Pi^2-\eta_\Pi)\overline{X}+\varphi_3(t,V)((\eta_\Pi^2-\eta_\Pi)^2+2(\eta_\Pi^2-\eta_\Pi))\overline{L}^2
+\varphi_4(t,V)(\eta_\Pi^2-\eta_\Pi)\overline{L}\\

&&\!\!+\varphi_5(t,V)((\eta_\Pi^2-\eta_\Pi)^2+2(\eta_\Pi^2-\eta_\Pi))\overline{X}\;\overline{L}\Big]=0.$$
\end{array}
\end{equation}
where$\varpi_1$, $\varpi_2$, $\varpi_3$ and $\varpi_4$ are given by Eq.(\ref{24}), Eq.(\ref{25}), Eq.(\ref{26}) and Eq.(\ref{27}), respectively. By Eq.(\ref{22}) and Eq.(\ref{23}), we have
\begin{equation}\label{39}
\pi^{\ast}(t,V) =-\dfrac{2\varphi_1(t,V)\overline{X}+\varphi_2(t,V)+\varphi_5(t,V)\overline{L}}{2\varphi_1(t,V)\overline{X}}\cdot\dfrac{\varpi_1}{\varpi_4}
-\dfrac{\varphi_5(t,V)\overline{L}}{2\varphi_1(t,V)\overline{X}}\cdot\dfrac{\varpi_2}{\varpi_4}-\dfrac{\varpi_3}{\varpi_4},
\end{equation}
where $\pi^{\ast}$ denotes the optimal investment decision on the risky asset. Substituting $\pi^{\ast}$ into Eq.(\ref{38}), and rearranging terms by the order of $\overline{X}^2$, $\overline{L}^2$ and $\overline{X}\;\overline{L}$, we obtain the following bivariate polynomial function of $\overline{X}$ and $\overline{L}$:
\begin{equation}
\begin{array}{rcl}\label{40}
$$&&\!\!\!\!\!\!\!\!\!\!\bigg[\varphi_{1t}+a_1(t)\varphi_{1}
+\kappa(\delta-V)\varphi_{1V}+\dfrac{1}{2}\sigma_V^2V\varphi_{1VV}+\lambda_V(\varphi_1(t,V+\eta_{VV})-\varphi_1(t,V))\\

&&\!\!\!\!\!\!\!\!\!\!+\lambda_c(\varphi_1(t,V+\eta_{Vc})-\varphi_1(t,V))+\lambda\beta_2^2\bigg]\overline{X}^2\\

&&\!\!\!\!\!\!\!\!\!\!+\bigg[\varphi_{2t}+a_2(t)\varphi_{2}
+\kappa(\delta-V)\varphi_{2V}+\dfrac{1}{2}\sigma_V^2V\varphi_{2VV}+\lambda_V(\varphi_2(t,V+\eta_{VV})-\varphi_2(t,V))\\

&&\!\!\!\!\!\!\!\!\!\!+\lambda_c(\varphi_2(t,V+\eta_{Vc})-\varphi_2(t,V))
+2\lambda(\alpha_2\beta_2-\beta_2^2X_2^*)\bigg]\overline{X}\\

&&\!\!\!\!\!\!\!\!\!\!+\bigg[\varphi_{3t}+(a_3(t)+V)\varphi_{3}
+(\kappa(\delta-V)+2\sigma_V\rho_{LV}V)\varphi_{3V}+\dfrac{1}{2}\sigma_V^2V\varphi_{3VV}+\lambda_V(\varphi_3(t,V+\eta_{VV})\\

&&\!\!\!\!\!\!\!\!\!\!-\varphi_3(t,V))+\lambda_c(\varphi_3(t,V+\eta_{Vc})-\varphi_3(t,V))
+\lambda_c\varphi_3(t,V+\eta_{Vc})(\eta_{Lc}^2+2\eta_{Lc})-\dfrac{\varphi_5^2}{4\varphi_1}\\

&&\!\!\!\!\!\!\!\!\!\!\cdot\dfrac{(\varpi_1+\varpi_2)^2}{\varpi_4}+\xi\varphi_5\bigg]\overline{L}^2\\

&&\!\!\!\!\!\!\!\!\!\!+\bigg[\varphi_{4t}+a_4(t)\varphi_{4}\!+\!(\kappa(\delta-V)
\!+\!\sigma_V\rho_{LV}V)\varphi_{4V}\!+\!\dfrac{1}{2}\sigma_V^2V\varphi_{4VV}\!+\!\lambda_V(\varphi_4(t,V+\eta_{VV})-\varphi_4(t,V))\\

&&\!\!\!\!\!\!\!\!\!\!+\lambda_c(\varphi_4(t,V+\eta_{Vc})-\varphi_4(t,V))+\lambda_c\varphi_4(t,V+\eta_{Vc})\eta_{Lc}
-\dfrac{\varphi_2\varphi_5}{2\varphi_1}\cdot\dfrac{\varpi_1(\varpi_1+\varpi_2)}{\varpi_4}+\xi\varphi_2\bigg]\overline{L}\\

&&\!\!\!\!\!\!\!\!\!\!+\bigg[\varphi_{5t}+a_5(t)\varphi_{5}+\lambda_c\eta_{Lc}\varphi_5(t,V+\eta_{Vc})
+(\kappa(\delta-V)+\sigma_V\rho_{LV}V)\varphi_{5V}+\dfrac{1}{2}\sigma_V^2V\varphi_{5VV}\\

&&\!\!\!\!\!\!\!\!\!\!+\lambda_V(\varphi_5(t,V+\eta_{VV})-\varphi_5(t,V))+\lambda_c(\varphi_5(t,V+\eta_{Vc})-\varphi_5(t,V))
+2\xi\varphi_1\bigg]\overline{X}\overline{L}\\

&&\!\!\!\!\!\!\!\!\!\!+\bigg[\varphi_{6t}-\lambda\varphi_{6}+\kappa(\delta-V)\varphi_{6V}+\dfrac{1}{2}\sigma_V^2V\varphi_{6VV}
+\lambda_V(\varphi_6(t,V+\eta_{VV})-\varphi_6(t,V))\\

&&\!\!\!\!\!\!\!\!\!\!+\lambda_c(\varphi_6(t,V+\eta_{Vc})-\varphi_6(t,V))
-\dfrac{\varphi_2^2}{4\varphi_1}\cdot\dfrac{\varpi_1^2}{\varpi_4}+\lambda(\alpha_2-\beta_2 X^*_2)^2\bigg]=0,$$
\end{array}
\end{equation}
where
\begin{equation}
\begin{array}{rcl}\label{41}
$$a_1(t)\!\!\!\!&=&\!\!\!\!\zeta^2-4\zeta\sigma_\Pi\rho_{\Pi r}-\lambda+2((m+\dfrac{\zeta^2}{2})-\mu_\Pi+\sigma_\Pi^2)
+\sigma_\Pi^2+\lambda_\Pi((\eta_\Pi^2-\eta_\Pi)^2\\

\!\!\!\!&{}&\!\!\!\!+2(\eta_\Pi^2-\eta_\Pi))-\dfrac{(\varpi_1+\varpi_3)^2}{\varpi_4},$$
\end{array}
\end{equation}
\begin{equation}
\begin{array}{rcl}\label{42}
$$a_2(t)\!\!\!\!&=&\!\!\!\!
(m+\dfrac{\zeta^2}{2})-\mu_\Pi+\sigma_\Pi^2-\zeta\sigma_\Pi\rho_{\Pi r}-\lambda+\lambda_\Pi(\eta_\Pi^2-\eta_\Pi)
-\dfrac{\varpi_1(\varpi_1+\varpi_3)}{\varpi_4},$$
\end{array}
\end{equation}
\begin{equation}
\begin{array}{rcl}\label{43}
$$a_3(t)\!\!\!\!&=&\!\!\!\!2(\mu_L-\mu_\Pi+\sigma_\Pi^2)+\sigma_{LS}^2+\sigma_\Pi^2-\lambda+\lambda_L\eta_L^2
\!+\!2\lambda_L\eta_L\!+\!\lambda_\Pi((\eta_\Pi^2-\eta_\Pi)^2+2(\eta_\Pi^2-\eta_\Pi)),$$
\end{array}
\end{equation}
\begin{equation}
\begin{array}{rcl}\label{44}
$$a_4(t)\!\!\!\!&=&\!\!\!\!\mu_L-\mu_\Pi+\sigma_\Pi^2-\lambda+\lambda_L\eta_L+\lambda_\Pi(\eta_\Pi^2-\eta_\Pi),$$
\end{array}
\end{equation}
\begin{equation}
\begin{array}{rcl}\label{45}
$$a_5(t)\!\!\!\!&=&\!\!\!\!(m+\dfrac{\zeta^2}{2})-2\mu_\Pi+3\sigma_\Pi^2-2\zeta\sigma_\Pi\rho_{\Pi r}
+\mu_L+\lambda_L\eta_L+\lambda_\Pi((\eta_\Pi^2-\eta_\Pi)^2+2(\eta_\Pi^2-\eta_\Pi))\\

\!\!\!\!&{}&\!\!\!\!-\lambda-\dfrac{\varpi_1^2+\varpi_1\varpi_2+\varpi_1\varpi_3+\varpi_2\varpi_3}{\varpi_4},$$
\end{array}
\end{equation}

Since Eq.(\ref{40}) holds for every $\overline{X}$ and $\overline{L}$, the following six PDEs holds with the boundary conditions:
\begin{equation}
\begin{array}{rcl}\label{46}
\left\{
\begin{aligned}
&\varphi_{1t}+a_1(t)\varphi_{1}
+\kappa(\delta-V)\varphi_{1V}+\dfrac{1}{2}\sigma_V^2V\varphi_{1VV}+\lambda_V(\varphi_1(t,V+\eta_{VV})-\varphi_1(t,V))&\\

&+\lambda_c(\varphi_1(t,V+\eta_{Vc})-\varphi_1(t,V))+\lambda\beta_2^2=0,&\\
&\varphi_1(T,V)=\beta_1^2,&
\end{aligned}
\right.
\end{array}
\end{equation}
\begin{equation}
\begin{array}{rcl}\label{47}
\left\{
\begin{aligned}
&\varphi_{2t}+a_2(t)\varphi_{2}
+\kappa(\delta-V)\varphi_{2V}+\dfrac{1}{2}\sigma_V^2V\varphi_{2VV}+\lambda_V(\varphi_2(t,V+\eta_{VV})-\varphi_2(t,V))&\\

&+\lambda_c(\varphi_2(t,V+\eta_{Vc})-\varphi_2(t,V))
+2\lambda(\alpha_2\beta_2-\beta_2^2X_2^*)=0,&\\

&\varphi_2(T,V)=2\beta_1(\alpha_1-\beta_1 X^*_1),&
\end{aligned}
\right.
\end{array}
\end{equation}
\begin{equation}
\begin{array}{rcl}\label{48}
\left\{
\begin{aligned}
&\varphi_{3t}+(a_3(t)+V)\varphi_{3}
+(\kappa(\delta-V)+2\sigma_V\rho_{LV}V)\varphi_{3V}+\dfrac{1}{2}\sigma_V^2V\varphi_{3VV}+\lambda_V(\varphi_3(t,V+\eta_{VV})&\\

&-\varphi_3(t,V))+\lambda_c(\varphi_3(t,V+\eta_{Vc})-\varphi_3(t,V))
+\lambda_c\varphi_3(t,V+\eta_{Vc})(\eta_{Lc}^2+2\eta_{Lc})-\dfrac{\varphi_5^2}{4\varphi_1}&\\

&\cdot\dfrac{(\varpi_1+\varpi_2)^2}{\varpi_4}
+\xi\varphi_5=0,&\\

&\varphi_3(T,V)=0,&
\end{aligned}
\right.
\end{array}
\end{equation}
\begin{equation}
\begin{array}{rcl}\label{49}
\left\{
\begin{aligned}
&\varphi_{4t}+a_4(t)\varphi_{4}+(\kappa(\delta-V)
+\sigma_V\rho_{LV}V)\varphi_{4V}+\dfrac{1}{2}\sigma_V^2V\varphi_{4VV}
+\lambda_V(\varphi_4(t,V+\eta_{VV})-\varphi_4(t,V))&\\

&+\lambda_c(\varphi_4(t,V+\eta_{Vc})-\varphi_4(t,V))+\lambda_c\varphi_4(t,V+\eta_{Vc})\eta_{Lc}
-\dfrac{\varphi_2\varphi_5}{2\varphi_1}\cdot\dfrac{\varpi_1(\varpi_1+\varpi_2)}{\varpi_4}+\xi\varphi_2=0,&\\

&\varphi_4(T,V)=0,&
\end{aligned}
\right.
\end{array}
\end{equation}
\begin{equation}
\begin{array}{rcl}\label{50}
\left\{
\begin{aligned}
&\varphi_{5t}+a_5(t)\varphi_{5}+\lambda_c\eta_{Lc}\varphi_5(t,V+\eta_{Vc})
+(\kappa(\delta-V)+\sigma_V\rho_{LV}V)\varphi_{5V}+\dfrac{1}{2}\sigma_V^2V\varphi_{5VV}&\\

&+\lambda_V(\varphi_5(t,V+\eta_{VV})-\varphi_5(t,V))+\lambda_c(\varphi_5(t,V+\eta_{Vc})-\varphi_5(t,V))
+2\xi\varphi_1=0,&\\

&\varphi_5(T,V)=0,&
\end{aligned}
\right.
\end{array}
\end{equation}
\begin{equation}
\begin{array}{rcl}\label{51}
\left\{
\begin{aligned}
&\varphi_{6t}-\lambda\varphi_{6}+\kappa(\delta-V)\varphi_{6V}+\dfrac{1}{2}\sigma_V^2V\varphi_{6VV}
+\lambda_V(\varphi_6(t,V+\eta_{VV})-\varphi_6(t,V))&\\

&+\lambda_c(\varphi_6(t,V+\eta_{Vc})-\varphi_6(t,V))
-\dfrac{\varphi_2^2}{4\varphi_1}\cdot\dfrac{\varpi_1^2}{\varpi_4}+\lambda(\alpha_2-\beta_2 X^*_2)^2=0,&\\
&\varphi_6(T,V)=(\alpha_1-\beta_1X^*_1)^2.&
\end{aligned}
\right.
\end{array}
\end{equation}

Next we solve the above equations from Eq.(\ref{46}) to Eq.(\ref{51}) one by one. First we solve Eq.(\ref{46}). Assume that $\widetilde{\varphi}_1(t,V)$ is the solution of the following system:
\begin{equation}
\begin{array}{rcl}\label{52}
\left\{
\begin{aligned}
&\widetilde{\varphi}_{1t}+a_1(t)\widetilde{\varphi}_{1}
+\kappa(\delta-V)\widetilde{\varphi}_{1V}+\dfrac{1}{2}\sigma_V^2V\widetilde{\varphi}_{1VV}
+\lambda_V(\widetilde{\varphi}_1(t,V+\eta_{VV})-\widetilde{\varphi}_1(t,V))&\\

&+\lambda_c(\widetilde{\varphi}_1(t,V+\eta_{Vc})-\widetilde{\varphi}_1(t,V))=0,&\\
&\widetilde{\varphi}_1(T,V)=\beta_1^2,&
\end{aligned}
\right.
\end{array}
\end{equation}
which has the following form
\begin{equation}\label{53}
\widetilde{\varphi}_1(t,V)=e^{\widetilde{\varphi}_{11}(t)+\widetilde{\varphi}_{12}(t)V},
\end{equation}
with terminal condition $\widetilde{\varphi}_1(T,V)=\beta_1^2$. Thus
\begin{equation}
\begin{array}{rcl}\label{54}
&&\widetilde{\varphi}_{1t}=(\widetilde{\varphi}_{11}'+\widetilde{\varphi}_{12}'V)\widetilde{\varphi}_{1},\qquad\qquad
\widetilde{\varphi}_{1V}=\widetilde{\varphi}_{12}\widetilde{\varphi}_{1},\qquad\qquad
\widetilde{\varphi}_{1VV}=\widetilde{\varphi}_{12}^2\widetilde{\varphi}_{1},\\
&&\widetilde{\varphi}_{1}(t,V+\eta_{VV})-\widetilde{\varphi}_{1}(t,V)=\big[e^{\widetilde{\varphi}_{12}(t)\eta_{VV}}-1\big]\widetilde{\varphi}_{1},\\
&&\widetilde{\varphi}_{1}(t,V+\eta_{Vc})-\widetilde{\varphi}_{1}(t,V)=\big[e^{\widetilde{\varphi}_{12}(t)\eta_{Vc}}-1\big]\widetilde{\varphi}_{1},\\
\end{array}
\end{equation}

Substituting Eq.(\ref{53}) and Eq.(\ref{54}) into Eq.(\ref{52}), we obtain
\begin{equation}\label{55}
\widetilde{\varphi}_{11}'+\widetilde{\varphi}_{12}'V+a_1(t)+\kappa(\delta-V)\widetilde{\varphi}_{12}
+\dfrac{1}{2}\sigma_V^2V\widetilde{\varphi}_{12}^2+\lambda_V\big[e^{\widetilde{\varphi}_{12}\eta_{VV}}-1\big]
+\lambda_c\big[e^{\widetilde{\varphi}_{12}\eta_{Vc}}-1\big]=0.
\end{equation}

Since the Eq.(\ref{55}) holds for every $V$, the following two equation systems hold:
\begin{equation}
\begin{array}{rcl}\label{56}
\left\{
\begin{aligned}
&\widetilde{\varphi}_{12}'-\kappa\widetilde{\varphi}_{12}+\dfrac{1}{2}\sigma_V^2\widetilde{\varphi}_{12}^2=0,&\\
&\widetilde{\varphi}_{12}(T)=0,&
\end{aligned}
\right.
\end{array}
\end{equation}
\begin{equation}
\begin{array}{rcl}\label{57}
\left\{
\begin{aligned}
&\widetilde{\varphi}_{11}'+a_1(t)+\kappa\delta\widetilde{\varphi}_{12}+\lambda_V\big[e^{\widetilde{\varphi}_{12}\eta_{VV}}-1\big]
+\lambda_c\big[e^{\widetilde{\varphi}_{12}\eta_{Vc}}-1\big]=0,&\\
&\widetilde{\varphi}_{11}(T)=\ln\beta_1^2.&
\end{aligned}
\right.
\end{array}
\end{equation}

Solving the above two systems, we have $\widetilde{\varphi}_{11}(t)=\ln\beta_1^2+\int_t^T a_1(s)ds$ and $\widetilde{\varphi}_{12}(t)=0$, thus, $\widetilde{\varphi}_1(t,V)$ is independent of the variable $V$, which can be written as
\begin{equation}\label{58}
\widetilde{\varphi}_1(t)=\beta_1^2 e^{\int_t^T a_1(s)ds},
\end{equation}
and the system (\ref{52}) can be rewritten as
\begin{equation}
\begin{array}{rcl}\label{59}
\left\{
\begin{aligned}
&\widetilde{\varphi}_{1t}+a_1(t)\widetilde{\varphi}_{1}=0,&\\
&\widetilde{\varphi}_1(T)=\beta_1^2.&
\end{aligned}
\right.
\end{array}
\end{equation}

Now we solve Eq.(\ref{46}). Let $T$ be a variable in $\widetilde{\varphi}_1$, i.e., $\widetilde{\varphi}_1(t)=\widetilde{\varphi}_1(t,T)=e^{\widetilde{\varphi}_{11}(t,T)}$, where $\widetilde{\varphi}_{11}(t,s)=\ln\beta_1^2+\int_t^s a_1(u)du$. We conjecture that
\begin{equation}\label{60}
\varphi_1(t)=\widetilde{\varphi}_1(t,T)+\bigg\{\displaystyle\int_t^T\widetilde{\varphi}_1(t,s)\lambda\beta_2^2ds\bigg\}\beta_1^{-2},
\end{equation}
thus
\begin{equation}\label{61}
\varphi_{1t}=\widetilde{\varphi}_{1t}+\Big[\displaystyle\int_t^T\widetilde{\varphi}_{1t}(t,s)\lambda\beta_2^2ds\Big]\beta_1^{-2}
-\lambda\beta_2^2.
\end{equation}

Substituting Eq.(\ref{60}) and Eq.(\ref{61}) into the left hand side of Eq.(\ref{46}), we have
\begin{equation}
\begin{array}{rcl}\label{62}
$$\!\!\!\!&{}&\!\!\!\!\widetilde{\varphi}_{1t}+\Big[\displaystyle\int_t^T\widetilde{\varphi}_{1t}(t,s)\lambda\beta_2^2ds\Big]\beta_1^{-2}
-\lambda\beta_2^2
+a_1(t)\bigg[\widetilde{\varphi}_1(t,T)+\Big[\displaystyle\int_t^T\widetilde{\varphi}_1(t,s)\lambda\beta_2^2ds\Big]\beta_1^{-2}\bigg]
+\lambda\beta_2^2\\

\!\!\!\!&=&\!\!\!\!\widetilde{\varphi}_{1t}+a_1(t)\widetilde{\varphi}_1(t,T)
+\Big[\displaystyle\int_t^T\widetilde{\big[\varphi}_{1t}(t,s)+a_1(t)\widetilde{\varphi}_1(t,s)\big]\lambda\beta_2^2ds\Big]\beta_1^{-2}\\

\!\!\!\!&=&\!\!\!\!0.$$
\end{array}
\end{equation}

Thus $\varphi_1(t)$ given by Eq.(\ref{60}) is the solution of system (\ref{46}), which is finally given by Eq.(\ref{30}). Similarly, $\varphi_2(t)$ is given by Eq.(\ref{31}) and $\varphi_6(t)$ is given by Eq.(\ref{35}).

Next we solve Eq.(\ref{50}). Since the coefficient of $\varphi_5$ and the constant term $2\xi\varphi_1$ are both independent of $V$, we suppose that $\varphi_5(t,V)$ is independent of $V$ and rewrite is as $\varphi_5(t)$. Thus Eq.(\ref{50}) can be rewritten as
\begin{equation}\label{63}
\varphi_{5t}+(a_5(t)+\lambda_c\eta_{Lc})\varphi_{5}+2\xi\varphi_1=0,
\end{equation}
with terminal value $\varphi_5(T)=0$, and $\varphi_5(t)$ is given by Eq.(\ref{34}).

Next we solve Eq.(\ref{49}). Since $\varphi_1$, $\varphi_2$ and $\varphi_5$ are all independent of $V$, and the coefficient $a_4$ is only dependent on time $t$, we assume that $\varphi_4$ is also independent of $V$, which satisfies the following system
\begin{equation}
\begin{array}{rcl}\label{64}
\left\{
\begin{aligned}
&\varphi_{4t}+(a_4(t)+\lambda_c\eta_{Lc})\varphi_{4}
-\dfrac{\varphi_2\varphi_5}{2\varphi_1}\cdot\dfrac{\varpi_1(\varpi_1+\varpi_2)}{\varpi_4}+\xi\varphi_2=0,&\\
&\varphi_4(T)=0,&
\end{aligned}
\right.
\end{array}
\end{equation}
thus $\varphi_4$ is given by Eq.(\ref{33}).

Since $\varphi_1$ and $\varphi_5$ are both independent of $V$, set
\begin{equation}\label{65}
f_3(t)=-\dfrac{\varphi_5^2}{4\varphi_1}\cdot\dfrac{(\varpi_1+\varpi_2)^2}{\varpi_4}+\xi\varphi_5.
\end{equation}

Let $\widetilde{\varphi}_{3}=\widetilde{\varphi}_{3}(t,V;\tau)$ be the solution of
\begin{equation}
\begin{array}{rcl}\label{66}
\left\{
\begin{aligned}
&\widetilde{\varphi}_{3t}(t,V;\tau)+(a_3(t)+V)\widetilde{\varphi}_{3}(t,V;\tau)
+(\kappa(\delta-V)+2\sigma_V\rho_{LV}V)\widetilde{\varphi}_{3V}(t,V;\tau)&\\

&+\dfrac{1}{2}\sigma_V^2V\widetilde{\varphi}_{3VV}(t,V;\tau)
+\lambda_V(\widetilde{\varphi}_3(t,V+\eta_{VV};\tau)-\widetilde{\varphi}_3(t,V;\tau))
+\lambda_c(\widetilde{\varphi}_3(t,V+\eta_{Vc};\tau)&\\

&-\widetilde{\varphi}_3(t,V;\tau))+\lambda_c\widetilde{\varphi}_3(t,V+\eta_{Vc};\tau)(\eta_{Lc}^2+2\eta_{Lc})=0,&\\

&\widetilde{\varphi}_3(\tau,V;\tau)=f_3(\tau),&
\end{aligned}
\right.
\end{array}
\end{equation}
and we have the following proposition.

\textbf{Proposition A} The solution of Eq.(\ref{48}) can be expressed as
\begin{equation}
\begin{array}{rcl}\label{67}
$$\varphi_3(t,V)=\displaystyle{\int_t^T}\widetilde{\varphi}_{3}(t,V;\tau)d\tau.$$
\end{array}
\end{equation}

\begin{proof} First we have $\varphi_3(T,V)=\displaystyle{\int_T^T}\widetilde{\varphi}_{3}d\tau=0.$ Set $\tau=t$ in the second equation of Eq.(\ref{66}), thus we have $\widetilde{\varphi}_3(t,V;t)=f_3(t)$. Differentiating Eq.(\ref{68}) with respect to $t$ and $V$, respectively, we have
\begin{equation}
\begin{array}{rcl}\label{68}
$$\varphi_{3t}\!\!\!\!&=&\!\!\!\!\displaystyle{\int_t^T}\widetilde{\varphi}_{3t}(t,V;\tau)d\tau-\widetilde{\varphi}_{3}(t,V;t)
=\displaystyle{\int_t^T}\widetilde{\varphi}_{3t}(t,V;\tau)d\tau-f_3(t),\\

\varphi_{3V}\!\!\!\!&=&\!\!\!\!\displaystyle{\int_t^T}\widetilde{\varphi}_{3V}(t,V;\tau)d\tau, \qquad
\varphi_{3VV}=\displaystyle{\int_t^T}\widetilde{\varphi}_{3VV}(t,V;\tau)d\tau.$$
\end{array}
\end{equation}

Substituting $\varphi_{3t}$, $\varphi_{3V}$ and $\varphi_{3VV}$ into Eq.(\ref{48})
\begin{equation}
\begin{array}{rcl}\label{69}
&&\displaystyle{\int_t^T}\widetilde{\varphi}_{3t}(t,V;\tau)d\tau-f_3(t)
+(a_3(t)+V)\displaystyle{\int_t^T}\widetilde{\varphi}_{3}(t,V;\tau)d\tau
+(\kappa(\delta-V)+2\sigma_V\rho_{LV}V)\\

&&\cdot\displaystyle{\int_t^T}\widetilde{\varphi}_{3V}(t,V;\tau)d\tau
+\dfrac{1}{2}\sigma_V^2V\displaystyle{\int_t^T}\widetilde{\varphi}_{3VV}(t,V;\tau)d\tau
+\lambda_V\displaystyle{\int_t^T}\widetilde{\varphi}_{3}(t,V+\eta_{VV};\tau)-\widetilde{\varphi}_{3}(t,V;\tau)d\tau\\

&&+\lambda_c\displaystyle{\int_t^T}\widetilde{\varphi}_{3}(t,V+\eta_{Vc};\tau)-\widetilde{\varphi}_{3}(t,V;\tau)d\tau
+\lambda_c\displaystyle{\int_t^T}\widetilde{\varphi}_{3}(t,V+\eta_{Vc};\tau)d\tau(\eta_{Lc}^2+2\eta_{Lc})
+f_3(t)\\

&&=\displaystyle{\int_t^T}0d\tau=0.
\end{array}
\end{equation}
\end{proof}

Now we start solving Eq.(\ref{66}). Suppose $\widetilde{\varphi}_{3}=\widetilde{\varphi}_{3}(t,V;\tau)=\widetilde{\varphi}_{31}(t;\tau)e^{\widetilde{\varphi}_{32}(t;\tau)V}$, with terminal value $\widetilde{\varphi}_{3}(\tau,V;\tau)=\widetilde{\varphi}_{31}(\tau)e^{\widetilde{\varphi}_{32}(\tau)V}=f_3(\tau)$. Thus
\begin{equation}
\begin{array}{rcl}\label{70}
$$&&\widetilde{\varphi}_{3t}=\bigg[\dfrac{\widetilde{\varphi}_{31}^{'}}{\widetilde{\varphi}_{31}}
+\widetilde{\varphi}_{32}^{'}V\bigg]\widetilde{\varphi}_{3},\qquad \widetilde{\varphi}_{3V}=\widetilde{\varphi}_{32}\widetilde{\varphi}_{3},\qquad \widetilde{\varphi}_{3VV}=\widetilde{\varphi}_{32}^2\widetilde{\varphi}_{3},\\

&&\widetilde{\varphi}_{3}(t,V+\eta_{VV};\tau)-\widetilde{\varphi}_{3}
=(e^{\widetilde{\varphi}_{32}\eta_{VV}}-1)\widetilde{\varphi}_{3},\qquad
\widetilde{\varphi}_{3}(t,V+\eta_{Vc};\tau)=e^{\widetilde{\varphi}_{32}\eta_{Vc}}\widetilde{\varphi}_{3}.$$
\end{array}
\end{equation}

Substituting Eq.(\ref{70}) into Eq.(\ref{66}) with the consideration of the terminal value, we obtain the following two systems
\begin{equation}
\begin{array}{rcl}\label{71}
\left\{
\begin{aligned}
&\widetilde{\varphi}_{32}^{'}+1+(2\sigma_V\rho_{LV}-\kappa)\widetilde{\varphi}_{32}+\dfrac{1}{2}\sigma_V^2\widetilde{\varphi}_{32}^2=0,&\\

&\widetilde{\varphi}_{32}(\tau)=0,&
\end{aligned}
\right.
\end{array}
\end{equation}
\begin{equation}
\begin{array}{rcl}\label{72}
\left\{
\begin{aligned}
&\dfrac{\widetilde{\varphi}_{31}^{'}}{\widetilde{\varphi}_{31}}+a_3+\kappa\delta\widetilde{\varphi}_{32}
+\lambda_V(e^{\widetilde{\varphi}_{32}\eta_{VV}}-1)+\lambda_c(e^{\widetilde{\varphi}_{32}\eta_{Vc}}-1)
+\lambda_c e^{\widetilde{\varphi}_{32}\eta_{Vc}}(\eta_{Lc}^2+2\eta_{Lc})=0,&\\

&\widetilde{\varphi}_{31}(\tau)=f_3(\tau).&
\end{aligned}
\right.
\end{array}
\end{equation}

We solve system (\ref{71}) first. Rewrite the first equation as
\begin{equation}\label{73}
\widetilde{\varphi}_{32}^{'}=-\dfrac{1}{2}\sigma_V^2\widetilde{\varphi}_{32}^2-(2\sigma_V\rho_{LV}-\kappa)\widetilde{\varphi}_{32}-1.
\end{equation}

Let $\triangle_3=(2\sigma_V\rho_{LV}-\kappa)^2-2\sigma_V^2$ be the discriminant of the following quadratic equation
\begin{equation}\label{74}
-\dfrac{1}{2}\sigma_V^2\widetilde{\varphi}_{32}^2-(2\sigma_V\rho_{LV}-\kappa)\widetilde{\varphi}_{32}-1=0.
\end{equation}

If $\triangle_3>0$, then the two real roots $h_{1,2}$ of Eq.(\ref{74}) can be expressed as
\begin{equation}\label{75}
h_{1,2}=\dfrac{(2\sigma_V\rho_{LV}-\kappa)-\sqrt{\triangle_3}}{\sigma_V^2}.
\end{equation}

Thus
\begin{equation}\label{76}
\widetilde{\varphi}_{32}(t)=\dfrac{h_1h_2e^{-\sqrt{\triangle_3}(\tau-t)}-h_1h_2}
{h_1e^{-\sqrt{\triangle_3}(\tau-t)}-h_2}.
\end{equation}

If $\triangle_3=0$, then we have
\begin{equation}\label{77}
\widetilde{\varphi}_{32}(t)=\dfrac{2\sigma_V\rho_{LV}-\kappa}{\sigma_V^2+\dfrac{1}{2}\sigma_V^2(\tau-t)(2\sigma_V\rho_{LV}-\kappa)}
-\dfrac{2\sigma_V\rho_{LV}-\kappa}{\sigma_V^2}.
\end{equation}

If $\triangle_3<0$, then
\begin{equation}\label{78}
\widetilde{\varphi}_{32}(t)=\sqrt{-\dfrac{\triangle_3}{\sigma_V^4}}\tan\Bigg[\arctan\bigg[\dfrac{2\sigma_V\rho_{LV}-\kappa}{\sqrt{-\triangle_3}}\bigg]
+\dfrac{1}{2}\sqrt{-\triangle_3}(\tau-t)\Bigg]-\dfrac{2\sigma_V\rho_{LV}-\kappa}{\sigma_V^2}.
\end{equation}

The solution of system (\ref{72}) is
\begin{equation}\label{79}
\widetilde{\varphi}_{31}(t)=e^{\int_t^{\tau} f_{31}(s)ds}\cdot f_3(\tau),
\end{equation}
where
\begin{equation}
\begin{array}{rcl}\label{80}
$$ f_{31}(t)=a_3+\kappa\delta\widetilde{\varphi}_{32}
+\lambda_V(e^{\widetilde{\varphi}_{32}\eta_{VV}}-1)+\lambda_c(e^{\widetilde{\varphi}_{32}\eta_{Vc}}-1)
+\lambda_c e^{\widetilde{\varphi}_{32}\eta_{Vc}}(\eta_{Lc}^2+2\eta_{Lc}),$$
\end{array}
\end{equation}
thus $\varphi_3(t,V)$ is given by Eq.(\ref{32}).

It is obvious that $2\varphi_1(t)>0$. Inserting $\varphi_1$, $\varphi_2$ and $\varphi_5$ into Eq.(\ref{39}), the optimal investment strategy is given by Theorem \ref{th3.1}.

\end{sloppypar}

\begin{thebibliography}{199}

\bibitem{Bjork}
Bj\"{o}rk, T.; Slinko, I. Towards a general theory of good-deal bounds. Review of Finance, \textbf{2006}, 10, 221--260.

\bibitem{Bodie}
Bodie, Z.; Detemple, J.B.; Otruba, S.; Walter, S. Optimal consumption-portfolio choices and retirement planning. Journal of Economic Dynamics and Control, \textbf{2004}, 28, 1115--1148.

\bibitem{Chen}
Chen, Z.; Li, Z.; Yan, Z.; Sun, J. Asset allocation under loss aversion and minimum performance
constraint in a DC pension plan with inflation risk. Insurance: Mathematics and Economics, \textbf{2017}, 75, 137--150.

\bibitem{Delong}
Delong, L.; Gerrard, R.; Haberman, S. Mean-variance optimization problems for an accumulation
 phase in a defined benefit plan. Insurance: Mathematics and Economics, \textbf{2008}, 42, 107--118.

\bibitem{Devolder} %% Books %% surname(s), initial(s), title, publisher, place of publication, year.
Devolder, P.; Janssen, J.; Manca, R.{\it Stochastic methods for pension funds}; Wiley:
New York, \textbf{2012}.

\bibitem{Dybvig}
Dybvig, P.H.; Liu, H. Lifetime consumption and investment: retirement and constrained borrowing. Journal of Economic Theory, \textbf{2010}, 145, 885--907.

\bibitem{Eisenberg}
Eisenberg, J. Optimal dividends under a stochastic interest rate. Insurance: Mathematics and Economics, \textbf{2015}, 65, 259--266.

\bibitem{Fleming} %% Books %% surname(s), initial(s), title, publisher, place of publication, year.
Fleming, W. H.; Soner, H. M. {\it Controlled Markov Processes and Viscosity Solutions}; Springer-Verlag:
New York, \textbf{1993}.

\bibitem{Guan}
Guan, G.; Liang, Z. Optimal management of DC pension plan in a stochastic interest rate and stochastic volatility framework. Insurance: Mathematics and Economics, \textbf{2014}, 57, 58--66.

\bibitem{Han}
Han, N.; Hung, M. Optimal asset allocation for DC pension plans under inflation. Insurance: Mathematics and Economics, \textbf{2012}, 51, 172--181.

\bibitem{XLiang}
Liang, X.; Bai, L.; Guo, J. Optimal time-consistent portfolio and contribution selection for
defined benefit pension schemes under mean-variance criterion. ANZIAM Journal, \textbf{2014}, 56, 66--90.

\bibitem{LiangandYuenandZhang16}Liang, Z.; Bi, J.; Yuen, K.C.; Zhang, C. Optimal mean--variance reinsurance and investment in a jump-diffusion financial market with common shock dependence. Mathematical Methods of Operations Research, \textbf{2016}, 84, 155--181

\bibitem{LiangandYuenandZhang18}
Liang, Z.; Yuen, K.C.; Zhang, C. Optimal reinsurance and investment in a jump-diffusion financial market with common shock dependence. Journal of Applied Mathematics and Computing, \textbf{2018}, 56, 637--664.

\bibitem{Merton}
Merton, R. C. Optimal consumption and portfolio rules in a continuous-time model. Journal of Economic Theory, \textbf{1971}, 3, 373--413.

\bibitem{Ngwira}
Ngwira, B.; Gerrard, R. Stochastic pension fund control in the presence of Poisson jumps. Insurance: Mathematics and Economics, \textbf{2007}, 40, 283--292.

\bibitem{Sun}
Sun, J.; Li, Z.; Zeng, Y. Precommitment and equilibrium investment strategies for defined
contribution pension plans under a jump-diffusion model. Insurance: Mathematics and Economics, \textbf{2016}, 67, 158--172.

\bibitem{Tang}
Tang, M.; Chen, S.; Lai, G. C.; Wu, T. Asset allocation for a DC pension fund under stochastic interest rates
and inflation-protected guarantee. Insurance: Mathematics and Economics, \textbf{2018}, 78, 87--104.

\bibitem{Wang}
Wang, S.; Lu, Y.; Sanders, B. Optimal investment strategies and intergenerational risk sharing for target benefit pension plans. Insurance: Mathematics and Economics, \textbf{2018}, 80, 1--14.

\bibitem{WangandLi}
Wang, P.; Lu, Z. Robust optimal investment strategy for an AAM of DC pension plans with stochastic rate and stochastic volatility. Insurance: Mathematics and Economics, \textbf{2018}, 80, 67--83.

\bibitem{Yao}
Yao, H.; Yang, Z.; Chen, P. Markowitz's mean-variance defined contribution pension fund
management under inflation: A continuous-time model. Insurance: Mathematics and Economics, \textbf{2013}, 53, 851--863.

\bibitem{Yong} %% Books %% surname(s), initial(s), title, publisher, place of publication, year.
Yong, X.; Zhou, X. Y. {\it Stochastic Controls: Hamiltonian Systems and {HJB} Equations}; Springer-Verlag:
New York, \textbf{1999}.

\bibitem{Zeng}
Zeng, Y.; Li, D.; Chen, Z.; Yang, Z. Ambiguity aversion and optimal derivative-based pension investment with stochastic income and volatility. Journal of Economic Dynamics and Control, \textbf{2018}, 88, 70--103.

\bibitem{Zhang07}
Zhang, A.; Ralf, K.; Ewald, C. Optimal management and inflation protection for defined contribution pension plans. Blatter der DGVFM, \textbf{2007}, 28, 239--258.

\bibitem{Zhang10}
Zhang, A.; Ewald, C. Optimal investment for a pension fund under inflation risk. Mathematical Methods of Operations Research, \textbf{2010}, 71, 353--369.

\bibitem{ZhangangGuo}
Zhang, X.; Guo, J. Optimal defined contribution pension management when risky asset and salary follow jump diffusion processes. East Asian Journal on Applied Mathematics, \textbf{2020}, 10, 22--39.

\end{thebibliography}
\end{document}